\documentclass[11pt]{article}
\usepackage{graphics}
\usepackage[demo]{graphicx}
\usepackage{epsfig}
\usepackage{pstricks}
\usepackage[normal]{subfigure}
\usepackage[latin1]{inputenc}
\usepackage[english,francais]{babel}
\usepackage{relsize,exscale}
\usepackage{makeidx}
\usepackage{enumitem}
\usepackage{amsfonts,amssymb,amsmath}
\usepackage{graphicx}
\usepackage{color}
\usepackage{multirow}
\usepackage{mathrsfs}
\usepackage[normalem]{ulem}
\usepackage{cancel}
\newenvironment{prooff}{{\it Proof :}}{\hfill\rule{2mm}{2mm}\vskip3mm\par}
\newtheorem{theorem}{Theorem}[section]
\newtheorem{lemma}[theorem]{Lemma}
\newtheorem{proposition}[theorem]{Proposition}
\newtheorem{corollary}[theorem]{Corollary}
\newtheorem{definition}[theorem]{Definition\rm}
\newtheorem{remark}{\it Remark\/}

\usepackage{color}
\definecolor{dred}{rgb}{0.92,0,0}
\definecolor{dgreen}{rgb}{0,0.92,0}
\definecolor{dblue}{rgb}{0,0,0.92}
\definecolor{dyellow}{rgb}{0.95,0.95,0}

\newcommand{\R}{\mathbb{R}}
\newcommand{\N}{\mathbb{N}}
\def\D{\displaystyle}
\newcommand{\hs}{\hspace{0.1cm}}

\newcommand{\sa}{\\ [0.2cm]}
\usepackage{multirow}
\graphicspath{
{./Figures/}
{./}
}
\title{On generalized binomial laws to evaluate finite element accuracy: \\
toward applications for adaptive mesh refinement}
\author{Jo\"el Chaskalovic \thanks{D'Alembert,
Sorbonne University, Paris, France, (\emph{email}: jch1826@gmail.com)}
\qquad
Franck Assous
\thanks{
Department of Mathematics, Ariel University, 40700 Ariel, Isra\"el, (\emph{email}: franckassous55@gmail.com).}
}
\date{}
\begin{document}
\maketitle
\selectlanguage{english}
\begin{abstract}
\noindent The aim of this paper is to provide new perspectives on the relative finite elements accuracy. Starting from a geometrical interpretation of the error estimate which can be deduced from Bramble-Hilbert lemma, we derive a probability law that evaluates the relative accuracy, considered as a random variable, between two finite elements $P_k$ and $P_m$, ($k < m$). We extend this probability law to get a cumulated probabilistic law for two main applications. The first one concerns a family of meshes and the second one is dedicated to a sequence of simplexes which constitute a given mesh. Both of this applications might be relevant for adaptive mesh refinement. \\[0.2cm]
\noindent{\footnotesize {\em keywords}: Error estimates, Probability, Finite elements, Bramble-Hilbert lemma, Mesh refinement.}
\end{abstract}
\section{Introduction}\label{intro}
\noindent The past decades have seen the development of finite element error estimates due to their influence on improving both accuracy and reliability in scientific computing. \sa
In the \emph{a priori} error estimates,  an unknown (in most of the cases) constant is involved which depends, among others, on the basis functions of the considered finite element and on a given semi-norm of the exact solution one wants to approximate. Moreover, error estimates are only upper bounds of the approximation error yielding that it exact value is unknown. \sa
This was the starting point  which motivated us to consider the approximation error as a random variable to therefore derive a probability law of the relative accuracy between two Lagrange finite elements $P_k$ and $P_m$, $(k<m)$, (\cite{Arxiv1} and \cite{Arxiv2}).\sa
Here, our aim is to extend the results obtained in \cite{Arxiv1} for a given mesh ${\mathcal T}_h$ defined by a fixed mesh size $h$ (corresponding to the largest diameter of all the simplexes of ${\mathcal T}_h$). We will generalize our purpose either to a family of meshes or, for a given mesh, to the sequence of the local elements that constitute it. \sa
As we will see, to proceed it, we will need to distinguish between the local and the global accuracy of two given Lagrange finite elements $P_k$ and $P_m$. \sa
Both of these points of view will be introduce as potential applications for adaptive mesh refinement.\sa
The paper is organized as follows. We recall in Section \ref{Geo_and_Proba} the mathematical problem we consider, the Bramble-Hilbert lemma and the resulting error estimate that allowed us to introduce a probability law for the relative global accuracy between $P_k$ and $P_m$ finite elements. Extension to a family of meshes is addressed in Section \ref{Family_meshes}. In Section \ref{from_global_to_local}, one proposes a generalization that describes the relative local accuracy between two finite elements in a sequence of simplexes belonging to a given mesh. Then, the cumulated probabilistic law is derived either for a family of meshes or for a sequence of simplexes. Concluding remarks follow.
\section{Comparing global finite elements accuracy by a probabilistic approach}\label{Geo_and_Proba}
\noindent Let $\Omega$ be an open bounded and non empty subset of $\R^{n}$ and $\Gamma$ its boundary which we assumed to be $C^1- $piecewise. Let also $u$ be the solution to the second order elliptic variational formulation:
\begin{equation}\label{VP}
\textbf{(VP}\textbf{)} \hspace{0.2cm} \left\{
\begin{array}{l}
\mbox{Find } u \in   V \mbox{ solution to:} \\ %[0.1cm]
a(u,v) = l(v), \quad\forall v \in V, %\\[0.2cm]
\end{array}
\right.
\end{equation}
where $V$ is a given Hilbert space endowed with the norm $ \left\|.\right\|_{V}$, $a(\cdot,\cdot)$ is a bilinear, continuous and $V-$elliptic form defined on $V \times V$, and $l(\cdot)$ a linear continuous form defined on~$V$. Classically, variational problem \textbf{(VP)} has one and only solution $u \in V$ (see for example \cite{ChaskaPDE}). In the sequel and for simplicity, we will restrict ourselves to the case where $V$ is a usual Sobolev space of distributions. \sa
Let us also consider the approximation $u_{h}$ of $u$, solution  to the approximate variational formulation:
\begin{equation}\label{VP_h}
\textbf{(VP}\textbf{)}_{h} \hspace{0.2cm} \left\{
\begin{array}{l}
\mbox{Find } u_{h} \in   V_h \mbox{ solution to:} \\ %[0.1cm]
a(u_{h},v_{h}) = l(v_{h}),\quad \forall v_{h} \in V_h, %\\[0.2cm]
\end{array}
\right.
\end{equation}
where $V_h$ is a given finite-dimensional subset of $V$. \sa
To state a corollary of Bramble-Hilbert's lemma (cf. \cite{Ciarlet} or \cite{RaTho82}), we assume that $\Omega$ is exactly covered by a mesh ${\mathcal T}_h$ composed by $N_K$ simplexes $K_{\mu}, (1 \leq \mu \leq N_K),$ which respect classical rules of regular discretization, (see for example \cite{ChaskaPDE} for the bidimensional case and \cite{RaTho82} in $\R^n$). Moreover, we denote by $P_k(K_{\mu})$ the space of polynomial functions defined on a given simplex $K_{\mu}$ of degree less than or equal to $k$, ($k \geq$ 1). \sa
Then, we have the following well-known result \cite{RaTho82}: \vspace{0.1 cm}
\begin{lemma}\label{Thm_error_estimate}
Suppose that there exists an integer $k \geq 1$ such that the approximation $u_h$ of $V_h$ is a continuous piecewise function composed by polynomials which belong to $P_k(K_{\mu}), (1\leq \mu\leq  N_K)$. \sa Then, $u_h$ converges to $u$ in $H^1(\Omega)$:
\begin{equation}
\D\lim_{h\rightarrow 0}\|u_h-u\|_{1,\Omega}=0.
\end{equation}
Moreover, if the exact solution $u$ belongs to $H^{k+1}(\Omega)$, we have the following error estimate:
\begin{equation}\label{estimation_error}
\|u_h-u\|_{1,\Omega} \hs \leq \hs \mathscr{C}_k\,h^k \, |u|_{k+1,\Omega}\,,
\end{equation}
where $\mathscr{C}_k$ is a positive constant independent of $h$, $\|.\|_{1,\Omega}$ the classical norm in $H^1(\Omega)$ and $|.|_{k+1,\Omega}$ denotes the semi-norm in $H^{k+1}(\Omega)$.
\end{lemma}
%\vspace{0.2cm}
%
In the sequel, we remind the probability law we derived in \cite{Arxiv1} which allowed us to evaluate the relative global accuracy, measured in $H^1$-norm, between two Lagrange finite elements.\sa
We consider two families of Lagrange finite elements $P_k$ and $P_m$ corresponding to a set of values $(k,m)\in \N^2$ such that $0 < k < m$. \\[0.1cm]
The two corresponding inequalities given by (\ref{estimation_error}), assuming that the solution $u$ to \textbf{(VP)} belongs to $H^{m+1}(\Omega)$, are respectively written as:
\begin{eqnarray}
\|u^{(k)}_h-u\|_{1,\Omega} \hs & \leq & \hs \mathscr{C}_k h^{k}\, |u|_{k+1,\Omega}, \label{Constante_01} \\%[0.1cm]
\|u^{(m)}_h\hspace{-0.09cm}-u\|_{1,\Omega} \hs & \leq & \hs \mathscr{C}_m h^{m}\, |u|_{m+1,\Omega}, \label{Constante_02}
\vspace{-1cm}
\end{eqnarray}
where $u^{(k)}_h$ and $u^{(m)}_h$ respectively denotes the $P_k$ and $P_m$ Lagrange finite element approximations of $u$, solution to $\textbf{(VP}\textbf{)}_{h}$.\\[0.2cm]
In what follows, for simplicity, we set $C_k\equiv \mathscr{C}_k \, |u|_{k+1,\Omega}$ and $C_m\equiv \mathscr{C}_m \, |u|_{m+1,\Omega}$. Therefore, inequalities (\ref{Constante_01}) and (\ref{Constante_02}) become:
\vspace{-0.2cm}
\begin{eqnarray}
\|u^{(k)}_h-u\|_{1,\Omega} \hs & \leq & \hs C_k h^{k}, \label{Constante_1} \\%[0.1cm]
\|u^{(m)}_h\hspace{-0.09cm}-u\|_{1,\Omega} \hs & \leq & \hs C_m h^{m} ,\label{Constante_2}
\end{eqnarray}
which show that the two polynomial curves defined by $f_k(h)\equiv C_k h^k$ and $f_m(h)\equiv C_m h^m$ are the upper bounds of the possible values for the two norms $\|u^{(k)}_h-u\|_{1,\Omega}$ and $\|u^{(m)}_h-u\|_{1,\Omega}$. More precisely, inequality (\ref{Constante_1}) (resp. (\ref{Constante_2})) indicates that the norm $\|u^{(k)}_h-u\|_{1,\Omega}$ (resp.  $\|u^{(m)}_h-u\|_{1,\Omega}$) is below the curve $f_k(h)$ (resp. $f_m(h)$).
\sa
As we are interested in comparing the relative positions of these two curves, we introduce their intersection point $h^*$ defined by:
\begin{equation}\label{h*}
\D h^* \equiv \left( \frac{C_k}{C_m}\right)^{\frac{1}{m-k}} = \left( \frac{\mathscr{C}_k \, |u|_{k+1,\Omega}}{\mathscr{C}_m \, |u|_{m+1,\Omega}}\right)^{\frac{1}{m-k}}.
\end{equation}
Now, as often in numerical analysis, there is no {\em a priori} information to surely or better specify the distance between $\|u^{(k)}_h-u\|_{1,\Omega}$ (resp. $\|u^{(m)}_h-u\|_{1,\Omega}$) and the curve $f_k$ or its precise value in the interval $[0, C_k h^k]$ due to (\ref{Constante_1}) (resp. the curve $f_m$ and the interval $[0, C_m h^m]$ due to (\ref{Constante_2})). \sa
Indeed, this situation is the consequence of two main ingredients:
\begin{enumerate}
\item The solution $u$ of \textbf{(VP)} is unknown, which motivates the use of a $P_k$ finite element method to build an approximation $u^{(k)}_h$,
\item The way the mesh generator processes the mesh is \emph{a priori} random which leads to a corresponding random approximation $u^{(k)}_h$ too.
\end{enumerate}
\noindent It is the reason why we treat the possible values of the norm $\|u^{(k)}_h-u\|_{1,\Omega}$ as a random variable defined as follows: \sa
Let us consider an experiment where the constitution of a random grid ${\mathcal T}_h$ and the corresponding random approximation $u^{(k)}_h$ are involved. Therefore, the approximation error $\|u^{(k)}_h-u\|_{1,\Omega}$ can also be viewed as a random variable, defined by the following probabilistic framework:
\begin{itemize}
\item A {\em random trial} corresponds to the grid constitution and the associated approximation $u^{(k)}_h$.
\item The probability space ${\bf\Omega}$ contains therefore all the possible results for a given random trial, namely, for all of the possible grids that the mesh generator may processed, or equivalently, for all the corresponding approximations $u^{(k)}_h$.
\end{itemize}
Then, for a fixed value of $k$, we define the random variable $X^{(k)}$ by:
\begin{eqnarray}
X^{(k)} : & {\bf\Omega} & \hspace{0.1cm}\rightarrow \hspace{0.2cm}[0,C_k h^k] \nonumber \\%[0.2cm]
& \boldsymbol{\omega}\equiv u^{(k)}_h & \hspace{0.1cm} \mapsto \hspace{0.2cm}\D X^{(k)}(\boldsymbol{\omega}) = X^{(k)}(u^{(k)}_h) = \|u^{(k)}_h-u\|_{1,\Omega}. \label{Def_Xi_h}
\end{eqnarray}
In the sequel, for simplicity, we will set: $X^{(k)}(u^{(k)}_h)\equiv X^{(k)}(h)$. \sa
Now, regarding the absence of information concerning the more likely or less likely values of the norm $\|u^{(k)}_h-u\|_{1,\Omega}$ in the interval $[0, C_k h^k]$, we will assume that the random variable $X^{(k)}(h)$ has a uniform distribution on the interval $[0, C_k h^k]$ in the following meaning:
\begin{equation}\label{Xk_uniform}
\forall\, (\alpha,\beta)\in\R_{+}^{2},\, 0 \leq \alpha < \beta \leq C_k h^k : Prob\left\{X^{(k)}(h) \in [\alpha,\beta]\right\}=\frac{\beta-\alpha}{C_k h^k}.
\end{equation}
\noindent Equation (\ref{Xk_uniform}) means that if one slides the interval $[\alpha,\beta]$ anywhere in $[0, C_k h^k]$, the probability of the event $\D\left\{X^{(k)}(h) \in [\alpha,\beta]\right\}$ does not depend on where the interval $[\alpha,\beta]$ is in $[0, C_k h^k]$; this is the property of uniformity of the random variable $X^{(k)}$. \sa
We are now able to evaluate the probability of the event
\begin{equation}\label{objectif}
\left\{\|u^{(m)}_h-u\|_{1,\Omega} \leq \|u^{(k)}_h-u\|_{1,\Omega}\right\} \equiv \left\{X^{(m)}(h) \leq X^{(k)}(h)\right\}
\end{equation}
to estimate the relative global accuracy between two finite elements of order $k$ and $m$, $(k<m)$.\sa
Let first of all define the relative global accuracy between two Lagrange finite elements $P_k$ and $P_m, (k<m),$ as follows:
\begin{definition}\label{Def_Global_Accuracy}
Let $P_k$ and $P_m$ ($k<m$) be two Lagrange finite elements. \sa
Then, we will say that "$P_m$ is \underline{globally} more accurate than $P_k$" if:
\begin{equation}\label{X(i)_h}
\|u^{(m)}_h-u\|_{1,\Omega}  \leq \|u^{(k)}_h-u\|_{1,\Omega}.
\end{equation}
\end{definition}
We will recall now the probabilistic law established in \cite{ChasAs18_1} or \cite{Arxiv2} to get an estimation on the relative global accuracy between two Lagrange finite elements $P_k$ and $P_m, (k<m)$, for a fixed mesh size $h$.
\begin{theorem}\label{The_nonlinear_law}
Let $u \in H^{k+1}(\Omega)$ be the solution to (\ref{VP}), and $u^{(i)}_h, (i=k \mbox{ or } i=m, k<m)$, the two corresponding Lagrange finite element $P_i$ approximations solution to (\ref{VP_h}).\\[0.1cm]
We assume the two corresponding random variables $X^{(i)}(h), (i=k \mbox{ or } i=m),$ defined by (\ref{Def_Xi_h}) are independent and uniformly distributed on $[0, C_i h^i]$, where $C_i$ are defined by (\ref{Constante_1})-(\ref{Constante_2}). \\[0.1cm]
Then, the probability such that "$P_m$ is \underline{globally} more accurate than $P_k$", as introduced in (\ref{Def_Global_Accuracy}), is given by:
\begin{equation}\label{Nonlinear_Prob}
\D Prob\left\{ X^{(m)}(h) \leq X^{(k)}(h)\right\} = \left |
\begin{array}{ll}
\D \hs 1 - \frac{1}{2}\!\left(\!\frac{\!\!h}{h^*}\!\right)^{\!\!m-k} & \mbox{ if } \hs 0 < h \leq h^*, \\[0.5cm]
\D \hs \frac{1}{2}\!\left(\!\frac{h^*}{\!\!h}\!\right)^{\!\!m-k} & \mbox{ if } \hs h \geq h^*.
\end{array}
\right.
\end{equation}
\end{theorem}
\begin{remark}

Immediate interpretations of the non linear probability law (\ref{Nonlinear_Prob}) are available:
\begin{enumerate}
\item $P_m$ finite element is not only \underline{asymptotically} more accurate than $P_k$ as $k<m$, when $h$ goes to 0, as usually considered and as a consequence of error estimate (\ref{estimation_error}). \sa
    Indeed, for all $h \leq h^*$, the probability of  "$P_m$ is more accurate than $P_k$" belongs to $[0.5, 1]$. It means that $P_m$ is \underline{more likely} accurate than $P_k$ for all of these values of $h$, and not only for arbitrarily small values of $h$.
\item On the contrary when $k<m$, the finite element $P_k$ is more likely accurate when $h>h^*$. This new point of view allows us to recommend that, for specific situations like for adaptive mesh refinement, this finite element would be more appropriated, as long as one would be able to detect that $h>h^*$.
\end{enumerate}
\end{remark}
The next section is devoted to a possible application of the probabilistic law (\ref{Nonlinear_Prob}) to a family of meshes, for example, in the process of mesh refinement.

\section{Extension to a family of meshes}\label{Family_meshes}
\noindent The aim of this section is to extend to a family of meshes the previous results we recalled for a given mesh ${\mathcal T}_h$. \sa
For this purpose, we introduce a family of $N$ regular meshes $\left({\mathcal T}^{(n)}\right)_{\!n=1,N}$ made of simplexes, each mesh ${\mathcal T}^{(n)}$ being characterized by its mesh size $h_n$. Let us also consider two Lagrange finite elements $P_k$ and $P_m$, $(k<m)$. \sa
Therefore, for each mesh ${\mathcal T}^{(n)}$, we can write the corresponding probability law (\ref{Nonlinear_Prob}) for the event $\D \left\{X^{(m)}(h_n) \leq X^{(k)}(h_n)\right\}$, associated to the given mesh size $h_n$. \sa
Our aim is now to evaluate the probability such that exactly $n_e$ meshes, $(n_e=0,1,\dots,N),$ satisfies "$P_m$ \emph{is more accurate than} $P_k$". To this end, let us introduce the sequence of $N$ independent Bernoulli random variables $(Y_n)_{n=1,N}$ defined by:
\begin{equation}\label{Bernoulli}
\D Y_n = \left |
\begin{array}{ll}
\hs 1 & \mbox{ if } X^{(m)}(h_n) \leq X^{(k)}(h_n), \medskip \\
\hs 0 & \mbox{ if not,}
\end{array}
\right.
\end{equation}
and also the random variable $S_N$ determined by:
\begin{equation}\label{SN_famille_maillages}
\D S_N = \sum_{n=1}^{N}Y_n.
\end{equation}
As each Bernoulli variable $Y_n$ indicates if "$P_m$ is more accurate than $P_k$" or not on the corresponding mesh ${\mathcal T}^{(n)}$, $S_N$ describes the number of meshes on the all $N$ meshes such that "$P_m$ is more accurate than $P_k$".
\begin{remark}
For any meshes ${\mathcal T}^{(p)}$ and ${\mathcal T}^{(q)}$ which belong to $({\mathcal T}^{(n)})_{n=1,N}$, characterized by their mesh size $h_p$ and $h_q$, the knowledge of the event "$P_m$ is more accurate than $P_k$" on ${\mathcal T}^{(p)}$ does not enable us to conclude anything on ${\mathcal T}^{(q)}$. \sa Hence, the $N$ random Bernoulli variables $Y_n, (1 \leq n \leq N)$, are considered as independent. \sa
It comes out from the nature of the event "$P_m$ is more accurate than $P_k$". Indeed,  for an unknown exact solution $u$ and the corresponding two approximations $u^{(i)}_{h_p}$ and $u^{(i)}_{h_q}, (i=k \mbox{ or } i=m),$ one cannot link the value $Y_p$ and $Y_q$, associated to a mesh size $h_p$ and $h_q$.
\end{remark}
Therefore, we have the following result:
\begin{proposition}\label{Distrib_exact_famille_maillages}
The distribution of probabilities corresponding to the \underline{exact} number of meshes satisfying "$P_m$ is more accurate than $P_k$" is given by:
\begin{eqnarray}
\hspace{-0.2cm} Prob\left\{S_N=0\right\}\,\,\, & = & \,\,\,\,(1-\mathcal{P}(h_1)) \dots (1-\mathcal{P}(h_N)), \label{F1} \\[0.2cm]
%\D \forall n=1,\dots,N-1 & : & \nonumber \\%[0.2cm]
\hspace{-0.2cm} Prob\left\{S_N=n_e\right\} & = & \!\!\!\!\!\!\!\!\!\!\!\!\!\!\!\sum_{ \substack{(i_1,\dots,i_N)\in \{1,\dots,N\} \\[0.1cm]
\mbox {\scriptsize and } \,i_l \neq i_q\, \mbox{\scriptsize for }  l \neq q
} }
\!\!\!\!\!\!\!\!\!\!\!\!\!\mathcal{P}(h_{i_1})..\mathcal{P}(h_{i_n})..(1-\mathcal{P}(h_{i_{n+1}}))..(1-\mathcal{P}(h_{i_N})), 1 \leq n_e \leq N-1, \label{F2}\\
\hspace{-0.2cm} Prob\left\{S_N=N\right\} & = & \,\,\,\,\,\mathcal{P}(h_1) \dots \mathcal{P}(h_N), \label{F3}
\end{eqnarray}
where the quantities $\mathcal{P}(h_{i_j})\equiv Prob\, \{ X^{(m)}(h_{i_j}) \leq X^{(k)}(h_{i_j})\}$ are given by the probability distribution (\ref{Nonlinear_Prob}) of theorem~\ref{The_nonlinear_law}.
\end{proposition}
Later in the paper, the sequence of values $\D \left(\frac{}{}\!\!Prob\left\{S_N=n_e\right\}\right)_{n_e=0,N}$ will be called the exact probabilistic distribution of the relative \emph{global accuracy}. \sa
\begin{prooff}
To establish formulas (\ref{F1})-(\ref{F3}), it is sufficient to notice that the $N$ random Bernoulli variables $(Y_n)_{n=1,N}$ defined by (\ref{Bernoulli}) are independent. \sa
Indeed, (\ref{F1})-(\ref{F3}) is a direct generalization of the binomial law for the variable $S_N$ we would have to consider if all the Bernoulli variables $(Y_n)_{n=1,N}$ had been defined by the same probability $p$, given by: $$\D p\equiv Prob \{Y_n=1\}=Prob \{X^{(m)}(h_n) \leq X^{(k)}(h_n)\}, \forall n=1 \mbox{ to } N.$$
\end{prooff}
\begin{remark}
From proposition \ref{Distrib_exact_famille_maillages} one can also get the \underline{cumulated} distribution of the number of meshes for which "$P_m$ \emph{is more accurate than} $P_k$". Namely it corresponds to the probability such that at least $n$ meshes, $(n=1,\dots,N)$, are such that "$P_m$ \emph{is more accurate than} $P_k$". \sa
Using the definition of the random variable $S_N$ given by (\ref{SN_famille_maillages}), it corresponds to the probabilities defined by:
\begin{equation}\label{Distrib_cumul_maillages}
\D \forall n_e=1,\dots,N: Prob\left\{S_N \geq n_e\right\} = \sum_{j=n_e}^{N}Prob\left\{S_N = j\right\}.
\end{equation}
\end{remark}
Formula (\ref{Distrib_cumul_maillages}) is not easy to explicitly expressed due to formulas (\ref{F1})-(\ref{F3}). However, in the next proposition, we will prove a recurrence relation which allows us to determine an algorithm to compute the exact probabilistic distribution as well as the cumulated one.\sa
To this end, let $p_N$ denotes the probability defined by $p_N \equiv Prob\left\{Y_N = 1\right\}$. Then, we have the following result:
\begin{proposition}\label{recurrence_exact_distrib}
Let $\left({\mathcal T}^{(n)}\right)_{\!n=1,N}$ be a family of $N$ regular meshes composed by simplexes, each mesh ${\mathcal T}^{(n)}$ being characterized by its mesh size $h_n$.\sa Then, we have:
\begin{equation}\label{recurrence_exact_distrib_Formula}
\forall\, n_e=1,N: Prob\left\{S_N = n_e\right\} \, = \, p_N\,Prob\left\{S_{N-1} = n_e-1\right\} + (1-p_N) \, Prob\left\{S_{N-1} = n_e\right\}.
\end{equation}
%
%where $p_N$ denotes the probability defined by: $p_N \equiv Prob\left\{Y_N = 1\right\}$.
%
\end{proposition}
\begin{prooff}
Formula (\ref{recurrence_exact_distrib_Formula}) corresponds to the decomposition of the event $\left\{S_N = n\right\}$ into two independent events which are:
$$\D\left(\!\!\frac{}{}\left\{S_{N-1} = n_e-1\right\} \cap \left\{Y_N = 1\right\}\right) \mbox{ and } \left(\!\!\frac{}{}\left\{S_{N-1} = n_e\right\} \cap \left\{Y_N = 0\right\}\right).$$
Then, (\ref{recurrence_exact_distrib_Formula}) results from classical probabilistic property.
\end{prooff}
\begin{remark}
Proposition \ref{recurrence_exact_distrib} enables us to process the computation of the exact probabilistic law which corresponds to $\D\left(\frac{}{}\!\!Prob\left\{S_N = n_e\right\}\right)_{n_e=0,N}$, and therefore, the cumulated one too.\sa
Indeed, first of all, let us notice that:
\begin{equation}\label{prob_1}
\D\forall n=1,N: Prob\left\{S_{n} = n\right\}=Prob\left\{\sum_{k=1}^{n}Y_k=1\right\}= \sum_{k=1}^{n}Prob\left\{Y_k=1\right\}= p_1\dots p_{n},
\end{equation}
as the bernoulli variables $(Y_k)_{k=1,n}, (\forall\, n=1,N),$ are independent.\sa
Moreover, with the same arguments, we also notice that:
\begin{equation}\label{prob_2}
\forall n=1,N: Prob\left\{S_{n} = 0\right\}=(1-p_1)\dots(1-p_{n}).
\end{equation}
Therefore, by the help of the recurrence relation (\ref{recurrence_exact_distrib_Formula}) and (\ref{prob_1})-(\ref{prob_2}), one can compute the exact probabilistic law $\D\left(\frac{}{}\!\!Prob\left\{S_N = n_e\right\}\right)_{n_e=0,N}$, step by step, as follows:\sa
Relations (\ref{prob_1})-(\ref{prob_2}) directly give: $\D Prob\left\{S_{n} = 0\right\}$ and $Prob\left\{S_{n} = n\right\}, \forall n=1,N$.
\begin{eqnarray}
\hspace{-0.8cm}\bullet \hspace{-0.2cm}&\mbox{Step 1} & : \hspace{0.2cm} Prob\left\{S_2 = 1\right\} \mbox{ by (\ref{recurrence_exact_distrib_Formula})},\nonumber\\[0.2cm]
\hspace{-0.8cm}\bullet \hspace{-0.2cm}&\mbox{Step 2} & : \hspace{0.2cm} Prob\left\{S_3 = 1\right\}, Prob\left\{S_3 = 2\right\} \mbox{ by (\ref{recurrence_exact_distrib_Formula})}, \nonumber\\[0.2cm]
\hspace{-0.8cm}& \dots & : \hspace{0.3cm}\dots, \hspace{0.3cm}\dots,\nonumber\\[0.2cm]
\hspace{-0.8cm}\bullet \hspace{-0.2cm}&\mbox{Step N} & : \hspace{0.2cm} Prob\left\{S_N = 1\right\}, \dots, Prob\left\{S_N = N-1\right\} \mbox{ by (\ref{recurrence_exact_distrib_Formula})},\nonumber
\end{eqnarray}
\end{remark}
\begin{remark}
The cumulated probabilistic distribution is then a direct consequence of proposition \ref{recurrence_exact_distrib} due to formula (\ref{Distrib_cumul_maillages}). Finally, in the particular case corresponding to $n=1$, we can derive an explicit expression of $Prob\left\{S_N \geq 1\right\}$. This is the purpose of the next proposition.
\end{remark}
\begin{proposition}\label{Prop_Couverture_mailages_N1_N2}
Let $N$ be the total number of meshes which belong to a family of regular meshes $\left({\mathcal T}^{(n)}\right)_{n=1,N}$ made of simplexes, and $h_n$ denotes the mesh size of a given mesh ${\mathcal T}^{(n)}$. We assume that $N = N_1 + N_2$, where $N_1$ is the number of meshes such that $h_n \leq h^*$ and $N_2$ those such that $h_n > h^* $. \sa
Then, we have:
\begin{equation}\label{Couverture_mailages_N1_N2}
\D Prob\left\{S_N \geq 1\right\} = 1- \frac{1}{2^{N_1}}\left[\frac{h_1 \dots h_{N_1}}{h^{*N_1}}\right]^{m-k}\left[1-\frac{1}{2}\left(\frac{h^*}{h_{N_1+1}}\right)^{m-k}\right]\dots\left[1-\frac{1}{2}\left(\frac{h^*}{h_{N}}\right)^{m-k}\right].
\end{equation}
\end{proposition}
\begin{prooff}
By definition of the opposite event of $\D \left\{S_N \geq 1\right\}$, we can directly write:
\begin{equation}\label{Couverture_maillages}
\D Prob\left\{S_N \geq 1\right\} = 1- Prob\left\{S_N = 0\right\} =  1 - (1-\mathcal{P}(h_1)) \dots (1-\mathcal{P}(h_N)),
\end{equation}
where we used the expression of (\ref{F1}). \sa
So, taking into account the probability law (\ref{Nonlinear_Prob}) of theorem \ref{The_nonlinear_law} we can explicit the probability such that on $N$ meshes, at least one mesh satisfies "$P_m$ \emph{is more accurate than} $P_k$". Indeed, we conclude by using the definitions of $N_1$ and $N_2$ and the corresponding expressions in the probabilistic law (\ref{Nonlinear_Prob}) to get (\ref{Couverture_mailages_N1_N2}).
\end{prooff}
\begin{remark}
As we recall before, when $k<m$, it is often sought that the $P_m$ finite element is more accurate than the $P_k$ one. However, the probability given by (\ref{Couverture_mailages_N1_N2}) shows that even the event "There is a least one mesh among $N$  meshes such that $P_m$ is more accurate than $P_k$" is not a sure event. Nevertheless, one can prove that this event is asymptotically sure. It is the purpose of the next proposition based on the following lemma.
\end{remark}
\begin{lemma}\label{Lemma_Convergence}
Let $\beta$ be a real number such that $0 < \beta < 1$. Let $p$ be a given integer and $(x_n)_{n\geq p}$ be a sequence of real numbers satisfying: $\D \forall\, n  \geq p, 0 < x_n \leq \beta$. \\ Then, the sequence product $(\Pi_N)_{N \in \N}$ defines by: $\D \Pi_N \equiv \prod_{n=p}^{N} x_n$ converges to~0 when $N$ goes to~$+\infty$.
\end{lemma}
\begin{prooff} $\frac{}{}$ \sa
$\blacktriangleright$ Convergence of $\Pi_N$. \\
As the sequence $(x_n)_{n\geq p}$ belongs to the interval $]0,\beta]$ such that $0 < \beta < 1$, we can write that:

\begin{equation}
\D 0 < \frac{\Pi_{N+1}}{\Pi_N} =\, x_{N+1} \leq \beta < 1.
\end{equation}
Therefore, the sequence $\Pi_N$ is a decreasing and bounded from below by 0, so it converges.\\[0.2cm]
$\blacktriangleright$  Limit of $\Pi_N$.\sa
To compute the limit of the sequence $\Pi_N$ let us consider the following inequalities:
\begin{equation}
\D 0 < x_p \dots x_N \leq \left(\max_{p \leq n \leq N} x_n \right)^{N-p+1} \leq \beta^{N-p+1} < 1,
\end{equation}
since we assume that $\beta$ belongs to $]0,1[$. \sa
As a consequence, the sequence $\left(\beta^{N-p+1}\right)_{N\in\N}$ goes to 0 when $N$ goes to infinity, and due to Squeeze theorem \cite{Steward}, the sequence $\Pi_N$ too.
\end{prooff}
\noindent We are now in position to state the following proposition.
\begin{proposition}\label{Couv_Asympt_Famille_Maillages}
Let $({\mathcal T}^{(n)})_{n=1,N}$ be a given family of regular meshes and $(h_n)_{n=1,N}$ the corresponding sequence of mesh sizes. We assume that there exists $h_{\mbox{\tiny max}}\in R_{+}^{*}$ such that:
\begin{equation}\label{maillage_fini}
\forall n>0: h_n\leq h_{\mbox{\tiny max}}\, \mbox{ such that }\, 0 < \mathcal{P}(h_{\mbox{\tiny max}}) < 1.
\end{equation}
Let also $S_N$ be the random variable introduced in (\ref{SN_famille_maillages}). \sa
Then, we have:
\begin{equation}\label{Form_Couv_Asympt_Famille_Maillages}
\D \lim_{N\rightarrow +\infty}Prob\left\{S_N \geq 1\right\} = \lim_{N\rightarrow +\infty}\sum_{n_e=1}^{N}Prob\left\{S_N = n_e\right\} = 1.
\end{equation}
\end{proposition}
\begin{prooff}
We consider the expression of $Prob\left\{S_N \geq 1\right\}$ given by (\ref{Couverture_maillages}) and we introduce the sequence $(x_n)_{n=1,N}$ defined by: $x_n \equiv 1-\mathcal{P}(h_n)$. \sa
Then, one can write $Prob\left\{S_N \geq 1\right\}$ as follows:
\begin{equation}
\D Prob\left\{S_N \geq 1\right\} = 1 - \prod_{n=1}^{N}x_n.
\end{equation}
Now, one can check from formula (\ref{Nonlinear_Prob}) that $$\forall\, n=1,N: x_n=1-\mathcal{P}(h_n) \in\, ]0,1[,$$
since $\mathcal{P}(h_n) \neq 0 \mbox{ and } \mathcal{P}(h_n)\neq 1, \forall\,h_n, 0 < h_n\leq h_{\mbox{\tiny max}}$, as we assume (\ref{maillage_fini}).\sa
Consequently, as the probabilistic law $\mathcal{P}(h)$ defined by (\ref{Nonlinear_Prob}) is decreasing on $\R_+$, we have:
\begin{equation}
\forall n=1,N: \mathcal{P}(h_n) \geq \mathcal{P}(h_{\mbox{\tiny max}}),
\end{equation}
and finally, it exists a real number $\beta\equiv 1- \mathcal{P}(h_{\mbox{\tiny max}})$, $(0 < \beta < 1)$, such that:
$$\forall\, n=1,N: x_n=1-\mathcal{P}(h_n) \leq \beta < 1.$$
Therefore, from lemma \ref{Lemma_Convergence} we obtain the claimed asymptotic behavior of $Prob\left\{S_N \geq 1\right\}$.
\end{prooff}
\begin{remark}
We notice that in proposition \ref{Couv_Asympt_Famille_Maillages} the assumption (\ref{maillage_fini}) corresponds to any concrete application where only finite mesh sizes are considered.
\end{remark}
The idea of proposition \ref{Couv_Asympt_Famille_Maillages} is that for a large number of meshes, the probability to get at least one mesh such that "$P_m$ \emph{is more accurate than} $P_k$" goes to 1 with the number of meshes. This property could be taken into account for the adaptive mesh refinement, where a family of meshes is built to handle the large variations of the approximate solution. \sa
In other words, this result indicates that for a large number of meshes, one has to consider that $P_m$ finite elements will surely be more accurate, at least on one mesh. However, it is not a sufficient indication to motivate the implementation of these finite elements which are more expensive than the $P_k$ finite elements when $k<m$. \sa
Furthermore, to apply this approach to adaptive mesh refinement, one needs to transfer the results of this section which concern \emph{global} behaviors of a family of meshes to get \emph{local} information to refine a given mesh, mainly depending on the local gradient of the approximated solution. This is the purpose of the next section.
\section{From a global to a local accuracy comparison of finite elements}\label{from_global_to_local}
\subsection{A local probability law to compare accuracy of finite elements}
\noindent In the previous sections, we described a probabilistic approach to estimate for a given mesh ${\mathcal T_h}$ or for a sequence of meshes $\left({\mathcal T}^{(n)}\right)_{n=1,N}$, the relative global accuracy between two Lagrange finite elements $P_k$ and $P_m, (k<m)$. \sa
We obtained laws of probabilities that depend on the mesh size $h$, namely the size of the largest diameter in the mesh. Accordingly,  the results we got are \emph{global} and do not explicitly consider any local behavior that could play an important role, particularly when one considers adaptive mesh refinement. \sa
Our purpose is now to derive a \emph{local} comparison tool between two Lagrange finite elements. This approach will be possible if one recalls that the error estimate (\ref{estimation_error}), deduced from Bramble Hibert's lemma, is elaborated by two main ingredients. The first one is Cea's lemma and the second one is the interpolation error (\emph{see} \cite{RaTho82}). \sa
So, for any simplex $K$ belonging to a regular mesh ${\mathcal T}_h$, let us introduce $\Pi^{(k)}_K$, the local $K-$ Lagrange interpolation operator of degree $k$, to define the local interpolation by the help of polynomials of degree lower or equal to $k$ on $K$.\sa
Then, one can write the global interpolation error $\|u-\Pi_h u\|_{1,\Omega}$ as follows:
\begin{equation}\label{Decomposition_Erreur_interpolation}
\D \|u-\Pi_h u\|_{1,\Omega} = \left(\sum_{K \in {\mathcal T}_h}\|u-\Pi^{(k)}_K u\|_{1,K}^2\right)^{1/2},
\end{equation}
where $\Pi_h$ denotes the global Lagrange interpolation operator on the mesh ${\mathcal T}_h$.\sa
So, under the conditions of theorem \ref{Thm_error_estimate}, for each $K \in {\mathcal T}_h$, to obtain the next local estimate, we follow P.A. Raviart and J. M. Thomas, (see \cite{RaTho82}), to get:
\begin{equation}\label{Estimation_erreur_locale_interpolation}
\|u-\Pi^{(k)}_K u\|_{1,K} \leq \mathscr{C}_k^{\,'} \,h_K^k \, |u|_{k+1,K}\leq \mathscr{C}_k^{\,'} \,h_K^k \, |u|_{k+1,\Omega}\,,
\end{equation}
where $h_K$ denotes the diameter of the simplex $K$ and $\mathscr{C}_k^{\,'}$ a positive constant which does not depend on $K$ and $h_K$, but depends on the reference element defining the Lagrange $P_k$ finite element \cite{RaTho82}. \sa
Then, the quantity $\mathscr{C}_k^{\,'}\,|u|_{k+1,\Omega}$ does not depend on $K$ too. This independency will be further crucial when we will extend our results to the mesh refinement process.\sa
Moreover, as a consequence of Cea's lemma, we can also deal with the following estimate:
\begin{equation}\label{Erreur_approx_Erreur_global_interpol}
\|u-u_h\|_{1,\Omega} \hs \leq \hs \frac{M}{\alpha}\,\|u-\Pi_h u\|_{1,\Omega},
\end{equation}
where $M$  is the continuity constant and $\alpha$ the ellipticity constant of the bilinear form $a(\cdot,\cdot)$. \sa
\noindent Now, due to (\ref{Decomposition_Erreur_interpolation}) and (\ref{Erreur_approx_Erreur_global_interpol}), we highlight that the accuracy of a given finite element $P_k$ can be \emph{locally} characterized by estimate (\ref{Estimation_erreur_locale_interpolation}). \sa
As a consequence, we define the relative \emph{local accuracy} between two finite elements $P_k$ and $P_m$ as the relative local interpolation accuracy on a given simplex $K$ as follows:
\begin{definition}\label{Def_Local_Accuracy}
Let $P_k$ and $P_m, (k<m),$ be two Lagrange finite elements and $K$ a given simplex which belongs to ${\mathcal T}_h$. We will say that "$P_m$ is \underline{locally} more accurate than $P_k$ on $K$" if:
\begin{equation}\label{X(i)_K_h}
\|u-\Pi^{(m)}_K u\|_{1,K} \leq \|u-\Pi^{(k)}_K u\|_{1,K}.
\end{equation}
\end{definition}
Therefore, if we assume that the exact solution $u$ of the variational formulation \textbf{(VP)} belongs to $H^{m+1}(\Omega)$, we can write inequality (\ref{Estimation_erreur_locale_interpolation}) for both of the $\Pi^{(i)}_K$, $(i=k \mbox{ or } i=m)$, local operators and we have:
\begin{eqnarray}
\|u-\Pi^{(k)}_K u\|_{1,K} \,\,& \leq & \mathscr{C}_k^{\,'} \,h_K^k \, |u|_{k+1,\Omega}\,, \label{Constante_k0} \\%[0.1cm]
\|u-\Pi^{(m)}_K u\|_{1,K} & \leq & \mathscr{C}_m^{\,'} \,h_K^m \, |u|_{m+1,\Omega}\,. \label{Constante_m0}
\vspace{-1cm}
\end{eqnarray}
So, if we set $C_k^{\,'}\equiv \mathscr{C}_k^{\,'} \, |u|_{k+1,\Omega}$ and $C_m^{\,'}\equiv \mathscr{C}_m^{\,'} \, |u|_{m+1,\Omega}$, inequalities (\ref{Constante_k0}) and (\ref{Constante_m0}) become:
\begin{eqnarray}
\|u-\Pi^{(k)}_K u\|_{1,K} \hs\hs & \leq & \hs C_k^{\,'} h_K^k, \label{Constante_k1} \\%[0.1cm]
\|u-\Pi^{(m)}_K u\|_{1,K} \hs & \leq & \hs C_m^{\,'} h_K^m, \label{Constante_m1}
\end{eqnarray}
which are the twins of inequalities (\ref{Constante_1}) and (\ref{Constante_2}). \sa
The main difference between (\ref{Constante_1})-(\ref{Constante_2}) and (\ref{Constante_k1})-(\ref{Constante_m1}) is the meaning of $h$. Here, in (\ref{Constante_k1})-(\ref{Constante_m1}), $h_K$ denotes the local diameter of the simplex $K$ whereas in (\ref{Constante_1})-(\ref{Constante_2}) $h$ is the involved maximum mesh size of ${\mathcal T}_h$. \sa
As a consequence, if we introduce the random variables $X^{(i)}_K(h_K), (i=k$ or $i=m$ and $k<m)$, defined by:
\begin{equation}\label{Def_Xi_K_h}
X^{(i)}_K(h_K) \equiv \|u-\Pi^{(i)}_K u\|_{1,K},
\end{equation}
we can directly get the probability of the event $\left\{X^{(m)}_K(h) \leq X^{(k)}_K(h)\right\}$ corresponding to "$P_m$ \emph{is \underline{locally} more accurate than} $P_k$ on $K$" defined in \ref{Def_Local_Accuracy} by adapting formulas (\ref{Nonlinear_Prob}) as follows:
\begin{corollary}\label{Prop_Local_Law_Proba}
Let $K$ be a given simplex of diameter $h_K$ belonging to a given regular mesh ${\mathcal T}_h$. Let $u\in H^{k+1}(\Omega)$ be the solution to (\ref{VP}) and $u^{(i)}_h, (i=k \mbox{ or } i=m, k<m)$, the corresponding Lagrange finite element $P_i$ approximations solution to (\ref{VP_h}).\\[0.1cm]
We assume the two corresponding random variables $X^{(i)}_K(h), (i=k \mbox{ or } i=m)$, defined by (\ref{Def_Xi_K_h}) are independent and uniformly distributed on $[0, C_i^{\,'} h_K^i]$, where $C_i^{\,'}$ are defined by (\ref{Constante_k1})-(\ref{Constante_m1}). \\[0.1cm]
Then, the probability such that "$P_m$ is \underline{locally} more accurate than $P_k$ on $K$" is given by:
\begin{equation}\label{Local_Nonlinear_Prob}
\D Prob\left\{ X^{(m)}_K(h) \leq X^{(k)}_K(h_K)\right\} = \left |
\begin{array}{ll}
\D \hs 1 - \frac{1}{2}\!\left(\frac{\!h_K}{h^*}\!\right)^{\!\!m-k} & \mbox{ if } \hs 0 < h_K \leq h^*, \medskip \\
\D \hs \frac{1}{2}\!\left(\frac{h^*}{\!h_K}\!\right)^{\!\!m-k} & \mbox{ if } \hs h_K \geq h^*,
\end{array}
\right.
\end{equation}
where $h^*$ is defined by (\ref{h*}), but where the constants $C_k^{\,'}$ and $C_m^{\,'}$ introduced in (\ref{Constante_k1}) and (\ref{Constante_m1}) replace $C_k$ and $C_m$ defined in (\ref{Constante_1}) and (\ref{Constante_2}).
\end{corollary}
\begin{remark}\label{Indep_K}
We notice that the corresponding value of $h^*$ does not depend on the simplex $K$ as we consider in inequality (\ref{Estimation_erreur_locale_interpolation}) the semi-norm of $u$ in $H^{k+1}(\Omega)$, on the one hand, and as the constant $\mathscr{C}_k^{\,'}$, (due to $\mathscr{C}_k$), does not depend on $K$ too, as we mention above, on the other hand.
\end{remark}
Therefore, formula (\ref{Local_Nonlinear_Prob}) gives us an evaluation of the \emph{local} accuracy comparison between two finite elements $P_k$ and $P_m$ based on the \emph{local} comparison accuracy between the corresponding $K-$ Lagrange interpolation errors of degrees $k$ and~$m$. \sa
We are now in position to extend this local result to a sequence of simplexes which belongs to a fixed mesh and which are concerned by adaptive mesh refinement process. This is the purpose of the next subsection.
\subsection{Toward applications for adaptive mesh refinement}
\noindent We now consider a \underline{given} mesh ${\mathcal T}_h$ composed by $N$ simplexes $K_\mu$ whose diameters are denoted by $(h_\mu)_{\mu=1,N}$.\sa
For each simplex $K_\mu, (1 \leq \mu \leq N),$ we consider the probability law of the event "$P_m$ \emph{is locally more accurate than} $P_k$ on $K_{\mu}$" which is given by (\ref{Local_Nonlinear_Prob}). Now, as in section \ref{Family_meshes}, let us introduce the $N$ Bernoulli random independent variables $(Y_\mu)_{\mu=1,N}$ defined by (\ref{Bernoulli}) where we replace $X^{(m)}(h)$ by $X^{(m)}_{K_{\mu}}(h)$, and also the corresponding random variable $S_N$ defined by (\ref{SN_famille_maillages}).\sa
Because the total similitude between the mathematical formalism of section \ref{Family_meshes} and the present one, we directly get the analogous formulas of (\ref{F1})-(\ref{F3}) and (\ref{Couverture_mailages_N1_N2}) but with a total different meaning. \sa
Particularly, by adapting (\ref{Couverture_mailages_N1_N2}) of proposition \ref{Prop_Couverture_mailages_N1_N2} to the current situation, we get the following result:
\begin{proposition}\label{Prop_Couverture_mailages_N1_N2_V2}
Let $N$ be the total number of simplexes of a given mesh ${\mathcal T}_h$ splitted so as $N=N_1+N_2$, where $N_1$ denotes the number of simplexes satisfying $h_\mu \leq h^*$ and $N_2$ such that $h_\mu > h^* $. Then, we have:
\begin{equation}\label{Couverture_mailages_N1_N2_V2}
\D Prob\left\{S_N \geq 1\right\} = 1- \frac{1}{2^{N_1}}\left[\frac{h_1 \dots h_{N_1}}{h^{*N_1}}\right]^{m-k}\left[1-\frac{1}{2}\left(\frac{h^*}{h_{N_1+1}}\right)^{m-k}\right]\dots\left[1-\frac{1}{2}\left(\frac{h^*}{h_{N}}\right)^{m-k}\right].
\end{equation}
\end{proposition}
\begin{remark}
We would like to highlight that even if the formalism is totally equivalent between this section and section \ref{Family_meshes}, one has to carefully distinguish the different meanings of the two situations. \sa Indeed, here the Bernoulli variables $(Y_\mu)_{\mu=1,N}$ determine for a given elementary mesh $K_\mu$ if "$P_m$ is \underline{locally} more accurate than $P_k$ on $K_\mu$", while the Bernoulli variables introduced in the previous section characterize if a given mesh belonging to a sequence of meshes is such that "$P_m$ is \underline{globally} more accurate than $P_k$" on this mesh.
\end{remark}
Proposition \ref{Prop_Couverture_mailages_N1_N2_V2} shows again that, albeit if it is usually assumed that finite elements $P_m$ are more accurate than $P_k$, ($k<m$), formula (\ref{Couverture_mailages_N1_N2_V2}) highlights that even "\emph{to get at least one simplex on} $N$ \emph{such that} $P_m$ \emph{is locally more accurate than} $P_k$" is not a sure event, as its probability is different from one. Moreover, one can get a similar result we got for proposition \ref{Couv_Asympt_Famille_Maillages}, but first of all, we have to adapt our framework to the adaptive mesh refinement process.\sa
So, let us now describe what happens for a given mesh in which adaptive mesh refinement is processed, to improve the solution in the areas such that the gradient of the approximated solution becomes large. \sa
For our purpose, we distinguish for a given mesh ${\mathcal T}_h$ two kinds of simplexes. Those who are not going to be changed, and those which will be refined. Let us denotes by $N', (N'<N),$ the number of new simplexes introduced in the mesh by the refinement process.\sa
We also introduce $N_1, (N_1\leq N'),$ the number of simplexes such that $h_\mu \leq h^*$ and $N_2=N'-N_1$ the rest of simplexes. Therefore, the random variable $S_N$ defined by (\ref{SN_famille_maillages}) becomes $S_{N'}$, which describes now the total number of new simplexes in the refinement process such that "$P_m$ is locally more accurate than $P_k$".\sa
We remark that these considerations make sense as we neutralize the dependency of $h^*$ toward any considered simplex $K$ which belongs to the given mesh ${\mathcal T}_h$, (see remark \ref{Indep_K}).\sa
Then, our interest is to determine the behavior of $\D Prob\left\{S_{N'} \geq 1\right\}$, (equivalently determined by (\ref{Couverture_mailages_N1_N2_V2}) if one changes $N$ by $N'$), when $N'$ goes to infinity, namely where the number of new simplexes $N'$ concerned by the refinement process becomes large. \sa
This situation corresponds to the framework of adaptive mesh refinement where one is usually interested by locally refine the mesh: this is performed by identifying in the mesh the areas such that the gradient of the approximated solution is high.\sa
Therefore, we are now in position to determine for a given set of new simplexes $N'$ which becomes large, the asymptotic probability such that "\emph{there exists at least one new simplex where} $P_m$ \emph{is locally more accurate than $P_k$}". \sa
This is the purpose of the following proposition which is the twin of proposition~\ref{Couv_Asympt_Famille_Maillages}.
\begin{proposition}\label{Couv_Asympt_Maillage_Adapt}
Let ${\mathcal T}_h$ be a given mesh composed by $N$ simplexes and let $N', (N'<N),$ be the number of new simplexes processed by a mesh refinement. We also assume that there exists $h_{\mbox{\tiny max}} \in \R^{*}_{+}$ which satisfies the equivalent condition of (\ref{maillage_fini}) for all the diameters $h_\mu, (\mu=1,N'),$ defining the new simplexes.\sa
Then we have:
\begin{equation}\label{Form_Couv_Asympt_Maillage_Adapt}
\D\lim_{N'\rightarrow +\infty}Prob\left\{S_{N'} \geq 1\right\} = 1.
\end{equation}
\end{proposition}
This result is not surprising due to the total similarity between formulas (\ref{Couverture_mailages_N1_N2}) and (\ref{Couverture_mailages_N1_N2_V2}), (if one remplaces $N$ by $N'$ in (\ref{Couverture_mailages_N1_N2_V2})). \sa
Again, it shows that one has to carefully consider the relative local accuracy between $P_m$ and $P_k$ Lagrange finite elements, $(k<m)$, except if locally, the number of simplexes becomes very large. \sa
However, for a fixed number $N'$ of new simplexes, this phenomena is more pronounced depending on the minimum number of simplexes $n'_e, (n'_e=1,N'),$ satisfying "$P_m$ would be more accurate than $P_k$". Unfortunately, deriving the analogous formulas of (\ref{F1})-(\ref{F3}) and (\ref{Distrib_cumul_maillages}) for a sequence of simplexes, to explicit the probability of the event $\{S_{N'}\geq n'_e\}$ is inextricable. \sa
So, one can compare for two different values of $n'_e$ the behavior of the corresponding probabilities. This is the purpose of the following proposition.
\begin{proposition}\label{SN'_n'1_n'2}
Let $n'_{e,1}$ and $n'_{e,2}$ be two integers such that $1 \leq n'_{e,1} < n'_{e,2} \leq N'$, and let $S_{N'}$ be the random variable equivalently defined as (\ref{SN_famille_maillages}).\sa
Then, we have:
\begin{equation}\label{SN'_n'1_n'2_Formula}
Prob\left\{S_{N'} \geq n'_{e,2}\right\} \leq Prob\left\{S_{N'} \geq n'_{e,1}\right\}.
\end{equation}
\end{proposition}
\begin{prooff}
The proof of proposition \ref{SN'_n'1_n'2} results from the following identity:
\begin{equation}
\D Prob\left\{S_{N'} \geq n'_{e,1}\right\} = \sum_{n'_e=\,n'_{e,1}}^{n'_{e,2}-1}Prob\left\{S_{N'} = n'_e \right\} + Prob\left\{S_{N'} \geq n'_{e,2}\right\},
\end{equation}
and consequently, (\ref{SN'_n'1_n'2_Formula}) holds.
\end{prooff}
Therefore, proposition \ref{SN'_n'1_n'2} clearly indicates that the more $n'_e$ the less the probability such that at least $n'_e$ simplexes statisfy "$P_m$ is locally more accurate than the $P_k$".\sa
This points out that we cannot easily disqualify the $P_k$ finite element in comparison with the $P_m$ one, particularly for a mesh refinement process.
\begin{remark}
We finally notice that it is also possible to compute the exact distribution and the associated cumulated one by applying the same principles we processed regarding the previous laws we derived for formula (\ref{recurrence_exact_distrib_Formula}) of proposition \ref{recurrence_exact_distrib} and for formulas (\ref{prob_1}) and (\ref{prob_2}) too.
\end{remark}
\section{Discussion and conclusion}\label{Conclusion}
\noindent In this paper, we proposed a new geometrical-probabilistic approach to evaluate the relative accuracy between two Lagrange finite elements $P_k$ and $P_m, (k<m)$. Basically, we distinguished two cases: a global approach and a local one. Both cases are based on a probabilistic interpretation of the error estimate one can derive from Bramble Hilbert lemma.\sa
In the global approach, we recall the probabilistic law (\ref{Nonlinear_Prob}) we derived in \cite{Arxiv1} to propose a first extension to a family of meshes. Regarding the local accuracy between two finite elements, we recall the two main components required to establish the \emph{a priori} error estimate and we highlighted that it is centrally based on the local interpolation error. This leads us to transpose our global analysis to the local one to get the corresponding probability distribution (\ref{Local_Nonlinear_Prob}).\sa
Afterwards, we extended our global and local results to two principal applications, both of them could concern with adaptive mesh refinement process. The global probability law has been used to describe the case of a family of meshes while the local probabilistic one helped us to consider a fixed mesh composed by a sequence of simplexes. \sa
These results strengthen those we got in \cite{Arxiv1} and show that even if we consider a family of meshes (respectively a sequence of simplexes), the event "\emph{to get at least one mesh (respectively, at least one simplex) such that} $P_m$ \emph{is more globally (respectively, locally) accurate than} $P_k$" is not a sure event, (cf. propositions \ref{Prop_Couverture_mailages_N1_N2}, \ref{Prop_Couverture_mailages_N1_N2_V2} and proposition \ref{SN'_n'1_n'2} for a more general case). \sa
\noindent However, it is proved in Propositions (\ref{Couv_Asympt_Famille_Maillages}) and (\ref{Couv_Asympt_Maillage_Adapt}) that, for a great number of meshes (respectively, a great number of simplexes), this event is asymptotically sure.\sa
As we noticed, these two applications - the family of meshes or the sequence of simplexes -  may be useful for adaptive mesh refinement.  We also proved a recurrence relation, (see formula (\ref{recurrence_exact_distrib})), in the case of a family of meshes, which can be adapted to the case of a sequence of simplexes. \sa
This enables to compute, for example, the probability such that "\emph{at least 50 percents of meshes (resp., at least 50 percents of simplexes) are such that} $P_m$ \emph{is more globally (resp., locally) accurate than} $P_k$". \sa
Finally, we have to mention that for all concrete applications, one will have to precisely estimate the critical value $h^*$. Indeed, all the probabilistic laws we derived are based on formulas (\ref{Nonlinear_Prob}) and (\ref{Local_Nonlinear_Prob}) which depend on $h^*$. \sa
Since $h^*$ strongly depends on the semi-norm $H^{k+1}(\Omega)$ of the exact solution $u$ to the variational problem, all the available techniques which belong to \emph{a priori} estimation theory of solutions to partial differential equations will be involved to evaluate $h^*$.\\[0.4cm]
\textbf{\underline{Acknowledgements}:}
The authors want to warmly dedicate this research to pay homage to the memory of Professors Andr\'e Avez and G\'erard Tronel who largely promote the passion of research and teaching in mathematics.
\end{document}